\documentclass[seceqn,number,dvips]{arxbj}

\aid{0}
\volume{16}
\issue{4}
\pubyear{2010}
\firstpage{1016}
\lastpage{1038}
\doi{10.3150/10-BEJ255}

\makeatletter

\newtheorem{lemma}{Lemma}[section]
\newtheorem{theorem}[lemma]{Theorem}

\newremark{remark}[lemma]{Remark}

\makeatother

\begin{document}
\begin{frontmatter}

\title{The limit distribution of the maximum increment of a
random walk with regularly varying jump size distribution}
\runtitle{The limit distribution of the maximum increment of a
random walk}

\begin{aug}
\author[a]{\fnms{Thomas} \snm{Mikosch}\thanksref{a}\ead[label=e1]{mikosch@math.ku.dk}} \and
\author[b]{\fnms{Alfredas} \snm{Ra\v ckauskas}\thanksref{b}\ead[label=e2]{alfredas.rackauskas@mif.vu.lt}\corref{}}
\runauthor{T. Mikosch and A. Ra\v ckauskas}
\pdfauthor{Thomas Mikosch, Alfredas Rackauskas}
\address[a]{Department of Mathematics,
University of Copenhagen, Universitetsparken 5, DK-2100
Copenhagen, Denmark. \printead{e1}}
\address[b]{Department of Mathematics, Vilnius University and Institute
of Mathematics and Informatics,\\ Naugarduko 24, LT-2006 Vilnius,
Lithuania. \printead{e2}}
\end{aug}

% HISTORY:
\received{\smonth{1} \syear{2009}}
\revised{\smonth{11} \syear{2009}}

% ABSTRACT
%
\begin{abstract}
In this paper, we deal with the asymptotic\ distribution\ of the
maximum increment
of a random walk with a regularly varying\ jump size distribution.
This problem
is motivated by a long-standing problem on change point detection
for epidemic alternatives. It turns out that the limit distribution\ of the
maximum increment of the random walk is one of the classical extreme
value distributions, the Fr\'echet distribution. We prove the results in the general
framework of point processes and for jump sizes taking values in a
separable Banach space.
\end{abstract}

% KEYWORDS
%

\pdfkeywords{Banach space valued random element, epidemic change
point, extreme value theory,
Frechet distribution, maximum increment of a random walk, point process
convergence, regular variation}

\end{frontmatter}

%s1 ###
\section{Introduction}\label{sec:1}
We commence by considering a sequence
$(X_i)$ of independent
%real-valued finite variance
random variables %with
%common distribution $F$
and denote the partial sums by
\[
S_0=0 ,\qquad S_n=X_1+\cdots+X_n ,\qquad n\ge1 .
\]
Our original goal is to investigate the asymptotic behavior of
the quantities
%
%e1.1 ###
\begin{equation}\label{eq:t-n}
T_n=\max_{1\leq\ell< n}\max_{0\le k\le n-\ell}
\bigl(\ell(1-\ell/n)\bigr)^{-1/2}
(S_{k+\ell}-S_k-\ell \overline X_n ) ,\qquad n\ge 1 ,
\end{equation}
where $\overline X_n$ denotes the sample mean of
$X_1,\ldots,X_n$. The normalization in $T_n$ is motivated by the
fact that, under the assumption of i.i.d. finite variance $X_i$,
$\operatorname{var}(S_{k+\ell}- S_k-\ell\overline X_n)$ is
proportional to $\ell
(1-\ell/n)$. In their book on change point analysis, Cs\"org\H{o}
and Horvath \cite{csorgohorvath1997} mention that nothing seems
to be known about the distributional properties of $T_n$.  There exist
several approaches to replace the original problem by a more
tractable one. One way is to restrict the range over which the
maximum is taken to {$\ell_n\leq\ell\leq n-\ell_n$} for some
$\ell_n\to\infty$ satisfying $\ell_n=\mathrm{o}(n)$; see, for example,
Yao \cite{yao1993}. Alternatively, one can change the normalizing
constants $\sqrt{\ell(1-\ell/n) }$ in a suitable way; see, for example,
Ra\v ckauskas and Suquet \cite{rackauskassuquet2004a}.

Statistics {of type} $T_n$ appear in the context of tests for
change points in the mean under epidemic alternatives. This
problem can be formulated as follows: given that $X_1, \dots,
X_n$ are independent random variables, test the null hypothesis of
constant mean
\begin{itemize}
\item $H_0\dvtx EX_1=EX_2=\cdots=EX_n=\mu$
\end{itemize}
against the \textit{epidemic alternative}
\begin{itemize}
\item $H_A\dvtx{}$There exist integers $1\le k^*<m^*<n$ such that
\begin{eqnarray*}
EX_1&=&\cdots=EX_{k^*}=EX_{m^*+1}=\cdots=EX_n=\mu,\\
EX_{k^*+1}&=&\cdots=EX_{m^*}=\nu\quad\mbox{and}\quad
\mu\ne\nu.
\end{eqnarray*}
\end{itemize}
One-sided alternatives such as $\mu>\nu$ or $\mu<\nu$ can also be
considered. Under the alternative $H_A$, the mean value
$\nu$ in the period $[k^\ast,m^\ast]$ is interpreted as an
epidemic deviation from the usual mean $\mu$ and $\ell^*=m^*-k^*$
is called the \textit{duration of the epidemic state.} To the best of
our knowledge, this kind of change point problem was formulated
for the first time by Levin and Kline \cite{levinkline1985} in
the context of abortion epidemiology. In the one-sided case, they
proposed the test statistic $\max_{1\le\ell\le n} \max_{0\le
k\le n-\ell} (S_{k+\ell}-S_k-\ell\overline X_n-\ell\delta/2)$, where
$\delta$ represents the smallest increment in the mean which is
sufficiently important to be detected. Simultaneously, epidemic-type models
were introduced by Commenges, Seal and Pinatel \cite{CSP1986}
in connection with experimental neurophysiology. They suggested a
circular representation of the model, allowing both $\ell^*$ and
$n-\ell^*$ to be interpreted as durations of the epidemic state.
Models with an epidemic-type change in the mean were also used for
detecting changed segments in non-coding DNA %deoxyribonucleic acid
sequences \cite{averyhenderson1999} and for studying
structural breaks in econometric contexts \cite{broemelingtsurumi}.% Moreover, the classical model
%of a box
%signal in additive noise, i.e.,
%X_i=a I_{[\theta_0, \theta_1]}(i/n)+\e_i,\ \ i=1, \ldots, n,
%for unknown parameters $a\in\bbr$ and
%$\theta_0, \theta_1\in(0,1)$ and iid mean zero noise $(\e_k)$
%also belongs to the class of epidemic-type models.
%{\bf Reference?}

The form of the test statistics $T_n$ is motivated by a
log-likelihood argument. Indeed, assuming $(X_i)$ to be i.i.d.
normal under
the hypothesis $H_0$ against the epidemic alternative $\mu<\nu$,
the test statistics $T_n$ is asymptotically equivalent to the square
root of a slightly generalized log-likelihood ratio statistics. In
the case of a two-sided epidemic alternative $\mu\ne\nu$, the
log-likelihood ratio statistics  under $H_0$ is asymptotically
equivalent to the quantity
%
%e1.2 ###
\begin{equation}\label{eq:twn} \widetilde T_n=\max_{1\le
\ell< n}\bigl(\ell(1-\ell/n)\bigr)^{-1/2} \max_{0\le k\le
n-\ell}|S_{k+\ell}-S_k-\ell\overline X_n| .
\end{equation}
Two-sided epidemic
alternatives, and hence test statistics such as $\widetilde T_n$, are
also meaningful in the { case} of multivariate observations
$X_i$. In this paper, we will even deal with sequences $(X_i)$ of
i.i.d. random elements with values in a separable Banach space.

Under the null hypothesis, when $\mu=EX_1$ is assumed to be known,
it is reasonable to replace the
sample mean $\overline X_n$ in the quantities $T_n$ and $\widetilde T_n$
by $\mu$. One then obtains the following
ramifications of $T_n$ and $\widetilde T_n$:
\begin{eqnarray}\label{eq:mwn}
{M}_n&=&\max_{1\le\ell\le n}\ell^{-1/2}
\max_{0\le k \le n-\ell} (S_{k+\ell}-S_k
-\ell\mu) ,\nonumber\\[-8pt]\\[-8pt]
\widetilde M_n&=&\max_{1\le\ell\le n}\ell^{-1/2}
\max_{0\le k\le n-\ell}
|S_{k+\ell}-S_k-\ell\mu| .\nonumber
\end{eqnarray}

Here, the choice of normalizing constants is again motivated by
the fact that the {variance $\operatorname{var}(S_{k+\ell}-S_k)$} is
proportional to $\ell$. An inspection of the quantities $T_n$,
$\widetilde T_n$, $M_n$ and $\widetilde M_n$ shows that under $H_0$,
{we may}
assume, without loss of generality, that the random variables $X_i$,
$i\ge1$,
have mean zero.

Various maximal elements of the random
field $(\ell^{-1/2}(S_{k+\ell}-S_k))_{\ell=1,\ldots,n,k=0,\ldots
,n-\ell}$
have
been widely discussed in the literature. Darling and Erd\H{o}s
\cite{darlingerdos1956} proved for a sequence\ $(X_i)$ of
i.i.d. standard normal random variables
and suitable constants $a_n>0$ and $b_n\in\mathbb{R}$ that
%
%e1.3 ###
\begin{equation}\label{eq:ed}
\lim_{n\to\infty} P \Bigl(a_n^{-1}\Bigl(\max_{\ell=1,\ldots,n}\ell
^{-1/2}S_\ell-b_n\Bigr)\le x \Bigr)=
\Lambda(x)=\mathrm{e}^{-\mathrm{e}^{-x}} ,\qquad x\in\mathbb{R}.
\end{equation}
Einmahl \cite{einmahl1989} showed that the
Darling--Erd\H{o}s result (\ref{eq:ed}) holds for suitable $a_n>0$
and $b_n\in\mathbb{R}$ if and only if $ E(X^2 I_{\{|X|\geq x\}
})=\mathrm{o}((\log\log
x)^{-1})$ as $x\to\infty$. The  limit distribution $\Lambda$ is
the \textit{Gumbel} or \textit{double exponential} extreme value distribution. { Note
that for a
sequence\ $(X_i)$ of i.i.d. standard normal random variables,} there
exist constants
$c_n>0$ and $d_n\in\mathbb{R}$ such that\
\[
\lim_{n\to\infty}P \Bigl(c_n^{-1}\Bigl(\max_{i=1,\ldots,n} X_i-d_n\Bigr)\le
x \Bigr)= \Lambda(x) ,\qquad x\in\mathbb{R};
\]
see \cite{gnedenko1943} and, for example, \cite{embrechtskluppelbergmikosch1997}, Example 3.3.29.
%Relation \eqref{eq:ed} means that, \textit{in an
%$(\ell^{-1/2}S_\ell)$ behave like the maxima of an iid \seq\ with \ds
%maximum domain of attraction of the Gumbel \ds.
A result in the same spirit was obtained by Siegmund and
Venkatraman \cite{siegmundvenkatraman1995}, who showed that for
i.i.d. standard normal random variables $X_i$, $i=1,2,\ldots,$ there exist
constants $a_n>0$ and $b_n\in\mathbb{R}$ such that\
$(a_n^{-1}(M_n-b_n))$ has a Gumbel limit distribution.
Another proof of this result is given in
\cite{kabluchko2008}. %{\bf give reference or cancel}
Finally, the famous Erd\H{o}s--R\'enyi laws are also closely
related to the maximum increments of a random walk. These laws
study the maxima of the random {sequence}
$(S_{k+\ell_n}-S_k)_{k=1,\ldots,n}$ for sequences $\ell_n\to
\infty$ with $\ell_n=\mathrm{o}(n)$; see, for example,
\cite{deheuvelsdevroye1987} for distributional convergence\ of the
Erd\H{o}s--R\'enyi statistic.\vspace*{1pt}
% and Arratia et al.
%in genetics.

In this paper, we are concerned with limit results for the
quantities $\widetilde T_n, T_n$ and $\widetilde M_n, M_n$, in the case where
$(X_i)$ is an i.i.d. sequence\ of heavy-tailed random variables. We
will obtain
results which parallel those in
\cite{siegmundvenkatraman1995} in the light-tailed case. A
useful definition of a heavy-tailed random variable\ $X$ with
distribution\ $F$ is
given via regular variation. The random variable\ $X$ is \textit{regularly
varying\ with index
$\alpha>0$} if there exists a slowly varying\ function\ $L$ such that\
$F$ satisfies the tail balance condition
%
%e1.4 ###
\begin{eqnarray}\label{eq:tb}
F(-x)\sim q \frac{L(x)}{x^\alpha} \quad\mbox{and}\quad
1-F(x)\sim p \frac{L(x)}{x^\alpha} ,\qquad x\to\infty,
\end{eqnarray}
where
$p\in(0,1)$, $p+q=1$; see \cite{binghamgoldieteugels1987} for an encyclopedic treatment
of regular variation.

Under the assumption of regular variation\ with index $\alpha$ on a generic
element $X$ of the i.i.d. sequence\ $(X_i)$, the theory developed in
Section \ref{sec:2} shows that the class of scaling factors
$\ell^{0.5}$ and $(\ell(1-\ell/n))^{0.5}$ which appear in the
quantities $\widetilde T_n, T_n$ and $\widetilde M_n, M_n$ is too narrow.
Indeed, the square root character of the normalizations suggests a
relationship with the central limit theorem, at least when
$\operatorname{var}(X)<\infty$. However, this argument is potentially
misleading. As a matter of fact, the scaling factors of the
maximum increments of such a random walk have to be chosen
depending on the index $\alpha$. They can range over a large class
of scaling functions. For $\gamma\ge0$, we define the class of functions
\begin{eqnarray*}
\mathcal{F}_\gamma&=&\{f\dvtx f \mbox{ is a positive non-decreasing
function
on }[0,\infty),  f(1)=1, f(\ell)\ge
\ell^\gamma,\\[-3pt]
&&\hphantom{\{}
\ell\ge1 \mbox{ and for any increasing sequence }(d_n)\mbox{ of positive
numbers}\\[-3pt]
&&\hphantom{\{}\mbox{such that }d_n^2/n\to0,\mbox{ it holds that
$\lim_{n\to\infty}\inf_{1\le\ell\le d_n}f(\ell
(1-\ell/n))/f(\ell)=1$}\} .
\end{eqnarray*}
Examples of functions in the class $\mathcal{F}_\gamma$ are
$f(x)=x^{\gamma'}$, where $\gamma'\ge\gamma$, and
$f(x)=x^\gamma\log^\beta(1+x),$ where $\beta>0$.

For any $f\in\mathcal{F}_\gamma$, we introduce the following
quantities:
\begin{eqnarray*}
\widetilde M_n^{(\gamma)}&=&\max_{1\le\ell\le n} (f(\ell))^{-1}
\max
_{0\le k\le n-\ell}
|S_{k+\ell}-S_k| ,\qquad
n\ge1 ,\\
\widetilde T_n^{(\gamma)}&=&\max_{1\le\ell< n}
\bigl(f \bigl(\ell(1-\ell/n) \bigr) \bigr)^{-1} \max_{0\le k\le
n-\ell}|S_{k+\ell}-S_k-\ell\overline X_n| ,\qquad
n\ge1 .
\end{eqnarray*}
We suppress the dependence of
the quantities $\widetilde M_n^{(\gamma)}$ and $\widetilde
T_n^{(\gamma)}$
on the function\ $f$. It will also turn out that the asymptotic\
results of this
section do not depend on the concrete form of the function\ $f$; they only
depend on the choice of $\gamma$.
We observe that $\widetilde M_n=\widetilde M_n^{(0.5)}$ for $f(\ell
)=\ell^{0.5}$ and
$\widetilde{T}_n=\widetilde{T}_n^{(0.5)}$ for $f(\ell)=\ell^{0.5}$
(cf. (\ref{eq:mwn}) and (\ref{eq:twn})).

The following result is a consequence\ of the general theory given in
Section \ref{sec:2}; see Theorem~\ref{thm:1}. In particular, the
result describes the asymptotic\ behavior of the quantities $\widetilde
M_n$ and
$\widetilde T_n$.
\begin{theorem}\label{thm:0} Consider an i.i.d. sequence\ $(X_i)$ of
random variables
which are regularly varying\ with index $\alpha>0$ and
have mean zero if it exists.
Then, for any function\ $f\in\mathcal{F}_\gamma$, $\gamma>\max
(0,0.5-\alpha^{-1})$,
%
%e1.5 ###
\begin{eqnarray}\label{eq:max}
\lim_{n\to\infty}P\bigl(a_n^{-1} \widetilde M_n^{(\gamma)}\le x\bigr) &=&
\Phi_\alpha(x)=\mathrm{e}^{-x^{-\alpha}} ,\qquad x>0 ,
\end{eqnarray}
where the
normalizing sequence\ is given by
%
%e1.6 ###
\begin{eqnarray}
\label{eq:an} a_n=\inf\{x\in
\mathbb{R}\dvtx P(|X|\leq x)\ge1-1/n\} .
\end{eqnarray}
Moreover,
%
%e1.7 ###
\begin{equation}\label{eq:max2}
\lim_{n\to\infty}P \bigl(a_n^{-1}
\widetilde T_n^{(\gamma)}\le x \bigr)= \Phi_\alpha(x) ,\qquad x>0 .
\end{equation}
\end{theorem}

Note that $\Phi_\alpha$ is the \textit{Fr\'echet}
extreme value distribution. In particular, for any i.i.d. sequence\ of
regularly varying\ random variables $X_i$
with index $\alpha>0$ and $(a_n)$ defined in (\ref{eq:an}),
\[
\lim_{n\to\infty}P\Bigl(a_n^{-1} \max_{i=1,\ldots,n} |X_i| \le x\Bigr)=
\Phi_\alpha(x) ,\qquad x\in\mathbb{R}.
\]
This relation follows
from classical results by Gnedenko \cite{gnedenko1943}; see, for example,
\cite{embrechtskluppelbergmikosch1997},
Theorem~3.3.7, for a more recent reference. Relation
(\ref{eq:max}) can thus be interpreted in the sense that the maximum of
the normalized increments $|S_{k+\ell}-S_k|$,
$\ell=1,\ldots,n$, $k=0,\ldots,n-\ell$, of the random walk
$(S_k)_{k=1,\ldots,n}$ is essentially determined by the maximum of
the i.i.d. random variables $|X_1|,\ldots,|X_n|$. The proof of
Theorem~\ref{thm:1}, in particular Lemma \ref{lem:2x}, explains
the asymptotic\ extreme value behavior. We mentioned above that Siegmund
and Venkatraman \cite{siegmundvenkatraman1995} proved an
analogous limit relation for $(a_n^{-1}(M_n-b_n))$, assuming that
$(X_i)$ is a
sequence\ of i.i.d. standard normal random variables. In this case,
the Gumbel distribution\ appears in the limit.
The scaling factor $\ell^{1-1/2}$ in $M_n$ is critical for
their result to hold. It can be interpreted as a boundary case
for the distributional limits of $(\widetilde M_n^{(\gamma)})$ when
$\alpha\to
\infty$. We also mention that the results in \cite{siegmundvenkatraman1995} go well beyond proving
convergence\ in distribution; they also give bounds for the probabilities
$P(a_n^{-1}(M_n-b_n)>x)$
as $n\to\infty$. Such bounds cannot be achieved with the methods used in
this paper.

It is worth mentioning that the statistics $\widetilde
T_n^{(\gamma)}$, with $f(x)=x^\gamma$ and $\gamma$ close to
$\max\{0, 1/2-1/\alpha\}$,
%$\alpha>2$,
allow one to detect epidemic changes in the mean, provided that the
duration of the epidemic state is of the order
$\ell^*=\mathrm{O}(n^\theta)$, where $\theta>\max\{1/\alpha,
2/(2+\alpha)\}$.
%for $1<\alpha\le2$, and $\theta>2/(2+\alpha)$ for $\alpha>2$.
Thus, for large $\alpha$, it is possible to detect short
epidemics. We refer to Cs\"org\H{o} and
Horvath \cite{csorgohorvath1997} and Ra\v ckauskas and
Suquet \cite{rackauskassuquet2004a} for details on applications
of statistics of the type $T_n$ to epidemic change problems.

The paper is organized as follows. In Section \ref{subsec:31}, we
introduce the notion of a regularly varying\ random element with values in
a Banach space and give several examples of such elements. The
main result of this paper (Theorem~\ref{thm:1}) is given in
Section \ref{subsec:32}. It proves that the normalized maximum
increment of a driftless random walk with values in a separable
Banach space and with regularly varying\ {jump sizes} converges in
distribution\ to
a Fr\'echet distribution. We complement this result with one-sided
versions for real-valued random variables. Section \ref{sec:3}
contains the
proofs of the results of Section \ref{sec:2}.
%Finally, in Section \ref{sec:4} we
%give some applications of these results to epidemic change point
%models.

%s2 ###
\section{General results}\label{sec:2}
In this section, we work in a framework more general than that of
Section \ref{sec:1}. Our generalizations are twofold:
(1) we consider i.i.d. sequences $(X_i)$ of Banach space valued,
regularly varying
\ random
elements; (2) we allow for more general normalizations of the
increments $S_{k+\ell}-S_k$, $\ell=1,\ldots, n$, $k=0,\ldots,n-\ell$.
In the following subsection, we introduce the notion of a regularly
varying\
random element and
in the subsequent subsection, we develop the asymptotic\ theory for
$\widetilde T_n$, $\widetilde M_n$ and related maximum increment quantities.

%s2.1 ###
\subsection{Regular variation in a Banach space}\label{subsec:31}
Consider %an iid \seq\ $(X_i)$ of random elements assuming
%values in
a separable Banach space $(\mathcal{B},\|\cdot\|).$ We say that a
$\mathcal{B}$-valued random element $X$ is \emph{regularly varying\
with index}
$\alpha>0$ if there exists a boundedly finite non-null measure\ $\mu$
on $ \mathcal{B}_0= \mathcal{B}\backslash\{\bf0\}$ such that
\[
\mu_n(\cdot)=n P(a_n^{-1} X \in\cdot)\stackrel{\widehat w}{\to} {
\mu(\cdot)} ,\qquad n\to\infty,
\]
where $\stackrel{\widehat w}{\to}$
is convergence\ in the sense that $\int_{\mathcal{B}_0}f\,\mathrm{d}\mu_n\to
\int_{\mathcal{B}_0}f\,\mathrm{d}\mu$ for any bounded and continuous function\ $f$ on $\mathcal
{B}_0$ with bounded support and where
%
%e2.1 ###
\begin{equation}\label{def:1}
a_n=\inf\{x\ge0\dvtx P(\|X\|\le x)\ge1-n^{-1}\}
\end{equation}
denotes the $(1-n^{-1})$-quantile of the distribution\ function of $\|
X\|$.
For locally compact $\mathcal{B}$, in particular for $\mathcal
{B}=\mathbb{R}^d$ for
some $d\ge1$, $\widehat w$-convergence\ coincides with vague
convergence\ and the
boundedly finite measures are the Radon measures; see
\cite{daleyverejones1988}, Appendix A2.6. The measure
$\mu$ necessarily satisfies the relation $\mu(t\cdot)= t^{-\alpha}
\mu(\cdot)$, $t>0$.
%which fact clarifies the relation with the index
%$\alpha$.
Moreover,
\[
P(\|X\|>x)=x^{-\alpha}L(x) \qquad\mbox{for a slowly varying\ function\ $L$} .
\]
In the case $\mathcal{B}=\mathbb{R}$, the notion of regular
variation\ of $X$
coincides with the definition given in (\ref{eq:tb}), provided that
$P(X>x)\sim p P(|X|>x)$ for some positive $p$. We refer to
Hult and Lindskog \cite{hultlindskog2006} for an insightful
survey on regular variation\ in complete separable metric spaces.
There, one
also finds a useful relation in terms of spherical coordinates
which is equivalent to regular variation\ of $X$ with index $\alpha
>0$: for
every $t>0$,
%
%e2.2 ###
\begin{equation}\label{eq:sph}
n P(\|X\|>t a_n ,X/\|X\|\in
\cdot)\stackrel{w}{\longrightarrow}t^{-\alpha} \widetilde P(\cdot),
\qquad n\to\infty,
\end{equation}
where $(a_n)$ is given by (\ref{def:1}) and $\widetilde P(\cdot)$
is a probability measure\ on the unit sphere $\mathbb{S}=\{x\in
\mathcal{B}\dvtx \|x\|=1\}$, called the \textit{spectral measure\ of $X$}, and $\stackrel
{w}{\longrightarrow}$
denotes weak convergence\ on the Borel $\sigma$-field $\mathcal
{B}_{\mathbb{S}}$ of
$\mathbb{S}$.

Examples of regularly varying\ random elements with values in a separable
Banach space can be found in, for example,
\cite{davismikosch2008} or \cite{ledouxtalagrand1993}. Those examples include max-stable random
fields on $[0,1]^d$ with a.s. continuous sample paths and regularly
varying\
finite-dimensional distributions \cite{davismikosch2008}. In this
case, the index $\alpha$
can be any positive number. Infinite variance stable processes
with values in a separable Banach space constitute another class of
regularly varying\ random elements; see
\cite{ledouxtalagrand1993}, Chapter 5, in particular Corollary 5.5.
In this
case, $\alpha$ is necessarily smaller than 2.

Another important example which is of interest in the context of
epidemic change point detection is a regularly varying sample
covariance operator, which we define next.
%The following example of \regvary\ random elements will be studied in
%Section \ref{sec:4}.
Denote the dual of $\mathcal{B}$ by $\mathcal{B}^*$ and let
$L(\mathcal{B}^*, \mathcal{B})$ be
the Banach space of {bounded} linear operators $u\dvtx\mathcal{B}^*\to
\mathcal{B}$ with norm
\[
\|u\|=\sup_{x^*\in\mathcal{B}^*\dvtx \|x^*\|\le1}
\|u(x^*)\| .
\]
For $x, y\in\mathcal{B}$, the operator $x\otimes y\dvtx
\mathcal{B}^*\to\mathcal{B}$ is defined by $(x\otimes y)
(x^*)=x^*(x)y$, $x^*\in
\mathcal{B}^*.$ It is immediate that $x\otimes y\in L(\mathcal
{B}^*,\mathcal{B})$ and
$\|x\otimes y\|=\|x\| \|y\|.$

Let $X$ be a $\mathcal{B}$-valued random element with mean zero and
finite second moment. The covariance
operator $\operatorname{cov}(X)=Q$ of $X$ then maps $\mathcal{B}^*$
into $\mathcal{B}$ and is
defined by
\[
Qx^*=E(x^*(X)X) ,\qquad x^*\in\mathcal{B}^* ,
\]
where the expectation is defined in the Bochner sense.

Assume that $X$ is defined on the probability space $(\Omega,
\mathcal{F}, P).$ Then, for each $\omega\in\Omega$, $(X\otimes
X)(\omega)=X(\omega)\otimes X(\omega)\in L(\mathcal{B}^*, \mathcal{B}).$
Evidently, $X\otimes X\dvtx \Omega\to L(\mathcal{B}^*, \mathcal{B})$ is
measurable,
that is, $X\otimes X$ is a random element with values in the
separable Banach space $L(\mathcal{B}^*, \mathcal{B}).$ To see this,
let $(X_n)$
be a sequence of { simple} functions that converge to $X$ a.s.
One then checks that $\|X_n\otimes X_n-X\otimes X\|\le
(\|X_n\|+\|X\|)\|X_n-X\|$ and, hence, $X_n\otimes X_n$ converge to
$X\otimes X$ a.s. Moreover, $X_n\otimes X_n$ are {simple}
functions in $L(\mathcal{B}^*, \mathcal{B}).$ For this random linear operator,
one can define regular variation\ with index $\alpha>0$ in the usual way.
%dicuss one possibility how one can achieve measurability of
%$X\otimes X$.}
We give the following result without a proof. It follows by means of a
standard continuous mapping argument.

\begin{lemma}\label{lem:7} If a
random element $X$ with values in $\mathcal{B}$ is
regularly varying\ with index $\alpha>0$, %and $X\otimes X$ is
%measurable
then $X\otimes X$ is regularly varying\ with index $\alpha/2.$
\end{lemma}

%s2.2 ###
\subsection{Results on the maximum increment of random
walks}\label{subsec:32}
%s2.2.1 ###
\subsubsection{Formulation of the main result}
Throughout this section, we consider an i.i.d. sequence\ of random
elements $X_i$, $i=1,2,\ldots,$ with values in a separable Banach
space $\mathcal{B}$. We assume that a generic element $X$ of this
sequence\ is regularly varying\ with index $\alpha>0$. If $\alpha>1$,
$E\|X\|<\infty$ and then its expectation $\mu=EX$ exists in the
Bochner sense. Since we are interested in quantities of the type
$\widetilde T_n$ and $\widetilde M_n$ defined in~(\ref{eq:twn}) and
(\ref{eq:mwn}), respectively, we {assume}, without loss of
generality, that $\mu=0$ whenever $\mu$ exists. Recall the
definition of the class of functions $\mathcal{F}_\gamma$, $\gamma
\ge0$,
from Section \ref{sec:1} and the definitions of the quantities
$\widetilde M_n^{(\gamma)}$ and $ \widetilde T_n^{(\gamma)}$ which we
adjust to
the case of Banach space valued random elements. Of course,
$f(\ell)=\ell^\gamma$, $\ell\ge0$, is a possible choice for
$f\in\mathcal{F}_\gamma$.

The following theorem is the main
result of this paper.

\begin{theorem}\label{thm:1} Let $(X_i)$ be a sequence of i.i.d.  %
random elements with values in a separable Banach space $\mathcal{B}$ and
assume that $X$ is regularly varying\ with index $\alpha>0$. In addition,
assume that $EX=0$ if $E\|X\|<\infty$ and
%
%e2.3 ###
\begin{equation}\label{restr:1}
\sup_{n\ge1}E\|n^{-1/\beta}S_n\|<\infty\qquad\mbox{for }
\cases{ \beta=2 , &\quad if $\alpha>2$,\cr
\mbox{every}  \beta<\alpha, &\quad if $1<\alpha\le 2$.
}
\end{equation}
Then, for $f\in\mathcal{F}_\gamma$,
$\gamma> \max(0, 0.5-\alpha^{-1})$, with the
normalizing sequence\ $(a_n)$ defined as in~(\ref{def:1}),
%
%e2.5 ###
%e2.4 ###
\begin{eqnarray}\label{main:1} \lim_{n\to\infty
}P\bigl(a_n^{-1}\widetilde
M_n^{(\gamma)}\le x\bigr)&= &\Phi_\alpha(x) ,\qquad x>0 ,\\
\label{main:2} \lim_{n\to\infty}P\bigl(a_n^{-1}\widetilde
T_n^{(\gamma)}\le x\bigr)&=& \Phi_\alpha(x) ,\qquad x>0 .
\end{eqnarray}
\end{theorem}

\begin{remark}\label{rem:1a}
It follows from (\ref{main:1}) and classical extreme value theory\ for
i.i.d. sequences that
the limit distributions of $(a_n^{-1} M_{n}^{(\gamma)})$ and
$(a_n^{-1} \max_{i=1,\ldots,n}\|X_i\|)$ coincide; see, for example,
\cite{embrechtskluppelbergmikosch1997},
Theorem~3.3.7. A theoretical explanation of this phenomenon is provided by
Lemma \ref{lem:2x}.
\end{remark}

The following result is the key to the proof of
Theorem~\ref{thm:1}. In particular, it explains why the quantities
\[
a_n^{-1}(f(\ell))^{-1} \max_{0\le k\le
n-\ell}\|S_{k+\ell}-S_k\| ,\qquad \ell\ge2 ,
\]
do not have any
influence on the limit behavior of $a_n^{-1}\widetilde M_n^{(\gamma)}$.
\begin{lemma}\label{lem:2x} Assume that $(X_n)$ is an i.i.d. sequence\ of
regularly varying\ random elements with values in a separable Banach space
$\mathcal{B}$ and index $\alpha>0$. The following statements
then hold:
\begin{enumerate}[(2)]
\item[(1)] for any $f\in\mathcal{F}_\gamma$, $\gamma\ge0$ and $h\ge1$,
\[
\lim_{n\to\infty}P \Bigl(a_n^{-1}\max_{1\le\ell\le h} (f(\ell))^{-1}
\max_{0\le
k\le n-\ell} \|S_{k+\ell}-S_k\| \le x \Bigr) = \Phi_\alpha(x) ,\qquad
x>0 ;
\]
\item[(2)] if we assume, in addition, that $EX=0$ if $E\|X\|<\infty$ and that
$(X_n)$ satisfies the condition~(\ref{restr:1}), then for any
$\delta>0$ and $f\in\mathcal{F}_\gamma$, $\gamma>\max(0,
0.5-\alpha
^{-1})$, we have
\[
\lim_{h\to\infty}\limsup_{n\to\infty}P \Bigl(\max
_{h\le\ell\le n}
(f(\ell))^{-1} \max_{0\le k\le n-\ell}\|S_{k+\ell}-S_k\|>\delta
a_n \Bigr)=0 .
\]
\end{enumerate}
\end{lemma}

%s2.2.2 ###
\subsubsection[A discussion of Theorem~2.2 and its
assumptions]{A discussion of Theorem~\protect\ref{thm:1} and its
assumptions}
In the following remarks, we provide a detailed discussion of the
statements and assumptions of Theorem~\ref{thm:1}.

\begin{remark}\label{rem:xu}
For $\gamma\ge1$, both relations (\ref{main:1}) and
(\ref{main:2}) are trivially satisfied. Indeed,
\[
\max_{0\le k\le n-1}\|X_k\|\le\widetilde M_n^{(\gamma)}\le\max
_{1\le
\ell\le n}(f(\ell))^{-1}\max_{0\le k\le
n-\ell}\sum_{i=k+1}^{k+\ell}\|X_k\|\le\max_{1\le k\le n}\|X_k\|
\]
for each $f\in\mathcal{F}_{\gamma}$ with $\gamma\ge1.$ (\ref
{main:1}) then follows in view of Remark \ref{rem:1a}.
\end{remark}

\begin{remark}\label{rem:97}
Under the assumptions of
Theorem~\ref{thm:1}, the sequence\ $(a_n^{-1}\widetilde M_n^{(\gamma
)})$ has
the same limit distribution\ as the sequence\
\[
a_n^{-1}\zeta_n^{(\gamma)}=a_n^{-1}\max_{1\le\ell\le
n}(f(\ell))^{-1} \max_{k=0,\ldots,n-\ell}\|S_{k+\ell}-S_k-\ell
\overline
X_n\| ,\qquad n\ge1 .
\]
To prove this statement, first
assume that $\gamma\ge1$. The argument of Remark \ref{rem:xu} then
shows that it suffices to consider the asymptotic\ behavior of the
sequence\
$(a_n^{-1}\max_{1\le k\le n}\|X_k-\overline X_n\|)$. If $E\|X\|<\infty
$, then
the strong law of large numbers\ ensures that\vspace*{-2pt} $\overline X_n\stackrel
{\mathrm{a.s.}}{\longrightarrow}EX$ as $n\to\infty$ and therefore
$a_n^{-1} \overline X_n\stackrel{P}{\longrightarrow}0$ as $n\to
\infty$. If $\alpha\in(0,1)$, then\vspace*{-2pt}
$a_n^{-1} \|\overline X_n\|\le a_n^{-1} n^{-1} \sum_{i=1}^n\|X_i\|
\stackrel{P}{\longrightarrow}
0$ as $n\to\infty$. In this case, $a_n^{-1}\sum_{i=1}^n \|X_i\|
\stackrel{d}{\longrightarrow}
Y_\alpha$ for some $\alpha$-stable random variable\ $Y_\alpha$ as
$n\to\infty$ since
$\|X\|$ is regularly varying\ with index $\alpha$; see
\cite{feller1971}, Section XVII, 5. The remaining case $\alpha=1$
with $E\|X\|=\infty$ is similar. In this case, again applying\vspace*{-2pt}
 \cite{feller1971}, Section XVII, 5, $a_n^{-1}\sum_{i=1}^n (\|X_i\|-
E(\|X\| I_{\{\|X\|\le a_n\}}))\stackrel{d}{\longrightarrow}Y_\alpha$ as
$n\to\infty$. Then, as $n\to\infty$,
\[
a_n^{-1} \|\overline X_n\|\le n^{-1}
a_n^{-1}\sum_{i=1}^n \bigl(\|X_i\|- E\bigl(\|X\| I_{\{\|X\|\le
a_n\}}\bigr)\bigr)+ a_n^{-1}E\bigl(\|X\| I_{\{\|X\|\le
a_n\}}\bigr)=\mathrm{o}_P(1) .
\]
In the last step, we also applied Karamata's theorem; see
\cite{binghamgoldieteugels1987}, Section 1.6.

We now consider the case $\gamma\in(0,1)$. Since $f\in\mathcal
{F}_\gamma$, it suffices to show that as $n\to\infty$,
%
%e2.6 ###
\begin{equation}\label{eq:67}
a_n^{-1} \max_{1\le\ell\le n} \ell^{-\gamma}\|\ell\overline X_n\|=
a_n^{-1} n^{1-\gamma}\|\overline X_n\|\stackrel{P}{\longrightarrow}0 .
\end{equation}
Assume that $\alpha>2$ and $\gamma>0.5-1/\alpha$. Then, by virtue of
(\ref{restr:1}), as $n\to\infty$,
\[
a_n^{-1} n^{1-\gamma}
\|\overline X_n\|=a_n^{-1} n^{-\gamma+0.5}\|n^{-0.5}S _n\|\stackrel
{P}{\longrightarrow}0 .
\]
If $\alpha\in(1,2]$, choose $\beta$ in $(\ref{restr:1})$
such that\ $\gamma>1/\beta-1/\alpha$. Then, as $n\to\infty$,
\[
a_n^{-1} n^{1-\gamma}
\|\overline X_n\|=a_n^{-1} n^{-\gamma+1/\beta}\|n^{-1/\beta}S _n\|
\stackrel{P}{\longrightarrow}0 .
\]
For $\alpha\in(0,1),$ we again use the fact that $(a_n^{-1}\sum
_{k=1}^n\|X_k\|)$
has an $\alpha$-stable limit as $n\to\infty$:
\[
a_n^{-1} n^{1-\gamma}
\|\overline X_n\|\le n^{-\gamma}a_n^{-1} \sum_{k=1}^n \|X_k\|
\stackrel{P}{\longrightarrow}0 .
\]
Similarly, for $\alpha=1$, as $n\to\infty$,
\[
a_n^{-1} n^{1-\gamma}
\|\overline X_n\|\le n^{-\gamma}
\Biggl[a_n^{-1} \Biggl(\sum_{k=1}^n \|X_k\| - n E\|X\|I_{\{\|X\|\le
a_n\}} \Biggr) +n^{1-\gamma} a_n^{-1} E\|X\|I_{\{\|X\|\le
a_n\}} \Biggr]
\stackrel{P}{\longrightarrow}0 .
\]
\end{remark}

\begin{remark}\label{rem:clt}
If $\alpha>2$, then condition (\ref{restr:1}) is fulfilled if the\vspace*{-2pt}
sequence\ $(X_i)$ satisfies the central limit theorem in $\mathcal{B}$,
that is, $n^{-1/2} S_n\stackrel{d}{\longrightarrow}Y$ as $n\to\infty$
for some Gaussian element in $\mathcal B$ (see Corollary 10.2 in
\cite{ledouxtalagrand1993}). If the space $\mathcal{B}$ is of type 2
(e.g., any finite-dimensional space, Hilbert space or Lebesgue
space $L_p$ with $p\ge2$), then (\ref{restr:1}) follows from regular
variation for $\alpha>2$. Similarly, if $\alpha\in(1,2)$ and $a_n^{-1}
S_n\stackrel{d}{\longrightarrow}Y_\alpha$ as $n\to \infty$ for some
$\alpha$-stable limit in $\mathcal{B}$, then (\ref{restr:1}) is
satisfied. This limit always exists in {a finite-dimensional space} as
a consequence\ of regular variation; see \cite{rvaceva1962}.
\end{remark}

\begin{remark}\label{rem:9}
The condition $\gamma>\max(0,
0.5-\alpha^{-1})$ divides the $\alpha$ values into two sets. For
$\alpha\le2$, this condition is satisfied for all $\gamma>0$,
whereas it restricts $\gamma$ to $(0.5-\alpha^{-1},\infty)$ for
$\alpha>2$. Under the assumption (\ref{restr:1}), this condition is
a natural one.  Indeed, assume for the moment that the $X_i$'s
are real-valued. By the definition of $(a_n)$,
$a_n=n^{1/\alpha}/\ell(n)$ for some slowly varying\ function\ $\ell$ and
hence
\[
a_n^{-1}\widetilde M_n^{(\gamma)}\ge
n^{-\alpha^{-1}-\gamma+0.5} \ell(n) |n^{-0.5}S_n| .
\]
If
$\gamma<0.5 -\alpha^{-1}$, then the left-hand side\ converges in
probability to
infinity {since $(n^{-0.5}S_n)$ converge in distribution to a
Gaussian random variable.} Hence, the normalization $(a_n)$ does
not ensure the stochastic boundedness of
$(a_n^{-1}\widetilde{M}_n^{(\gamma)}).$

With a
stronger normalization, a non-degenerate limit distribution of
the sequence\ $(\widetilde{M}_n^{(\gamma)})$ can be achieved
by an application of the invariance principle in H\"older space.
Following \cite{rackauskassuquet2006},
choose $f(\ell)=\ell^\gamma$ for $\gamma<0.5 -\alpha^{-1}$
and some $\alpha>2$
and assume that the central limit theorem for $(X_n)$ holds.
Then, as $n\to\infty$,
%
%e2.7 ###
\begin{equation}\label{eq:rwq}
n^{-0.5+\gamma}\widetilde{M}_n^{(\gamma)}\stackrel
{d}{\longrightarrow}R_{W,Q}=\sup_{s, t\in
[0, 1],
s\not=t}\frac{\|W_{Q}(t)-W_{Q}(s)\|}{|t-s|^\gamma} ,
\end{equation}
where $(W_{Q}(t))_{ 0\le t\le1}$ is the $\mathcal{B}$-valued Wiener process
corresponding to the covariance operator $Q=\operatorname{cov}(X).$
\end{remark}

\begin{remark}
Remark \ref{rem:9} shows that $\gamma=0.5-\alpha^{-1}$ for $\alpha
>2$ is the
borderline which divides the possible limit distributions of the normalized
sequence\
$(\widetilde M_n^{(\gamma)})$ into two classes: the Fr\'echet extreme
value distribution,\ as
described in
Theorem~\ref{thm:1}, and the distribution\ of the functional $R_{W,Q}$
of a
Wiener process given in (\ref{eq:rwq}). In the former case, only the
extremes in the sample $\|X_1\|,\ldots,\|X_n\|$ are responsible for
the
limit distribution, whereas in the latter case, the limit distribution\
is obtained by
an application of the
functional central limit theorem acting on the increments of the random
walk $(S_n)$.

The limit distribution\ of the normalized sequence\ $(\widetilde
M_n^{(0.5-\alpha^{-1})})$ is, in general, unknown; it very much depends
on the asymptotic\ behavior of the slowly varying\ function\ $L$ in the
tail $P(\|X\|>x)=x^{-\alpha}L(x)$. To illustrate the complexity of the
situation, we briefly consider two different cases. If $L(x)\to0$ as
$x\to\infty$, then the limiting relation (\ref{eq:rwq}) proved in
\cite{rackauskassuquet2006} still applies with $f(\ell)=\ell^\gamma$.
If $L(x)\sim c\in(0,\infty)$ as $x\to\infty$, then one can show that
$(a_n^{-1}\widetilde{M}_n^{(0.5-\alpha^{-1})})$ is stochastically
bounded. However, none of the sequences $(a_n^{-1}M_{n1})$ and
$(a_n^{-1} M_{n2})$ is asymptotically negligible in this case, where,
for any $h\ge1$ and $\varepsilon_n\to0$ as $n\to\infty$,
\begin{eqnarray*}
M_{n1}&=&\max_{1\le\ell\le h}\ell^{-0.5+\alpha^{-1}}\max_{0\le
k<n-\ell}\|S_{k+\ell}-S_k\| ,\\
M_{n2}&=&\max_{n\varepsilon_n\le\ell\le n}\ell^{-0.5+\alpha
^{-1}}\max_{0\le k<n-\ell}\|S_{k+\ell}-S_k\| .
\end{eqnarray*}
This is in stark contrast to the situation described in
Lemma \ref{lem:2x}. By the latter result, $(a_n^{-1}M_{n1})$
converges in distribution\ to a Fr\'echet-distributed random variable.
%%For the sake of
%illustration, we assume $\bbb=\bbr$.
On the other hand, by adapting\vspace*{-2pt}
the proof of Theorem~8 in \cite{rackauskassuquet2004a}, one can deduce that
$a_n^{-1}M_{n2}\stackrel{d}{\longrightarrow}R_{W,Q}$ as $n\to\infty
$ is possible, at least for $\mathcal{B}
=\mathbb{R}$.
\end{remark}

%s2.2.3 ###
\subsubsection{One-sided results for real-valued random variables}
In the remainder of this section, we restrict our attention to a
real-valued i.i.d. sequence\ $(X_i)$. We note that the two-sided
relations (\ref{eq:max}) and (\ref{eq:max2}) in
Theorem~\ref{thm:0} immediately follow from Theorem~\ref{thm:1}
by choosing $\mathcal{B}=\mathbb{R}$. However, in the real-valued
case, one can
also study one-sided versions of Theorem~\ref{thm:1}, for example, the
asymptotic\ behavior of the quantities, for $f\in\mathcal{F}_\gamma
$, $\gamma\ge0$ and $n\ge1$,
\begin{eqnarray*}
M_n^{(\gamma)}&=&\max_{1\le\ell\le n} (f(\ell))^{-1} \max_{0\le
k\le
n-\ell}(S_{k+\ell}-S_k) ,\\
m_n^{(\gamma)}&=&\min_{1\le\ell\le n} (f(\ell))^{-1} \min_{0\le
k\le
n-\ell}(S_{k+\ell}-S_k) .
\end{eqnarray*}

\begin{theorem}\label{thm:1a}
Assume that $(X_i)$ is an i.i.d. sequence\ of
real-valued random variables with distribution\ $F$ {which is regularly
varying} with index
$\alpha>0$, in the sense of (\ref{eq:tb}). In addition, assume that
$EX=0$ if the mean of $X$ exists.
%for $\alpha>1.$
Then, for $f\in\mathcal{F}_\gamma$, $\gamma>\max(0, 0.5-\alpha^{-1})$,
%
%e2.8 ###
\begin{eqnarray}\label{eq:limit}
&&\lim_{n\to\infty} P \bigl(({
p^{1/\alpha} a_n)^{-1}}m_n^{(\gamma)}\le-x ,
(p^{1/\alpha} a_n)^{-1} M_n^{(\gamma)}\le y
\bigr)\nonumber\\[-8pt]\\[-8pt]
&&\quad =
\Phi_\alpha(y)\bigl(1-\Phi_\alpha^{q/p}(x)\bigr) ,\qquad x,y>0 ,\nonumber
\end{eqnarray}
where the
normalizing sequence\ $(a_n)$ is defined as in (\ref{eq:an}) and $p\in
(0,1)$ appears in the tail balance condition (\ref{eq:tb}).
\end{theorem}
%
%defined as \beam\label{eq:0}
%b_n=\inf\{x\in\bbr: F(x)\ge1-1/n\} . \eeam\ethe\noindent
%{\bf Shall we give a result with centering by $\ell\ov X_n$?}
%
\begin{remark}
We note that as $n\to\infty$,
\[
(p^{1/\alpha} a_n)^{-1}\bigl(m_n^{(\gamma)},M_n^{(\gamma)}\bigr)\stackrel{d}{\longrightarrow}
\bigl(y^{(\gamma)},Y^{(\gamma)}\bigr) ,
\]
where the limit distribution is given
in (\ref{eq:limit}). In particular, $y^{(\gamma)}$ is independent
of $Y^{(\gamma)}$ and the range statistic
$M_n^{(\gamma)}-m_n^{(\gamma)}$ has the limit, as $n\to\infty$,
\[
(p^{1/\alpha} a_n)^{-1}\bigl(M_n^{(\gamma)}-m_n^{(\gamma)}\bigr)\stackrel
{d}{\longrightarrow}
Y^{(\gamma)}-y^{(\gamma)} .
\]
The limit distribution\ is the
convolution $\Phi_\alpha\ast\Phi_\alpha^{q/p}$, corresponding to
the sum of two independent Fr\'echet-distributed random variables.
\end{remark}

Consider the following one-sided version of the statistics
$\widetilde{T}_n^{(\gamma)}$:
\[
T_n^{(\gamma)}=\max_{1\le\ell< n} \bigl(\ell(1-\ell/n)\bigr)^{-\gamma}
\max_{0\le k\le n-\ell}(S_{k+\ell}-S_k-\ell\overline X_n) ,\qquad
n\ge1 .
\]

\begin{theorem}\label{thm:1b}
Assume that $(X_i)$ is an i.i.d. sequence\ of
real-valued random variables with distribution\ $F$ {which is regularly
varying} with index
$\alpha>0$, in the sense of (\ref{eq:tb}). In addition, assume that
$EX=0$ if the mean of $X$ exists. Then, for any $\gamma>\max(0,
0.5-\alpha^{-1})$,
%
%e2.9 ###
\begin{equation}\label{eq:limit1a}
\lim_{n\to\infty}P\bigl((p^{1/\alpha}a_n)^{-1}T_n^{(\gamma)}\le x\bigr)=
\Phi_\alpha(x) ,\qquad x>0 .
\end{equation}
\end{theorem}

The following quantity
has a structure similar to $M_n^{(\gamma)}$ for $f\in\mathcal
{F}_\gamma$:
\[
\widehat M_n^{(\gamma)}= \max_{\ell=1,\ldots
,n}(f(\ell))^{-1}
\max_{k=\ell+1,\ldots,n-\ell}(S_{k+\ell}+S_{k-\ell}-2 S_k) .
\]
In contrast to the quantities $M_n^{(\gamma)}$, where we
need to assume that $EX=0$ for $\alpha>1$ in order to guarantee the
asymptotic\ results of Theorem~\ref{thm:1a}, centering of the $X_i$'s in
$\widehat M_n^{(\gamma)}$ is automatic. Indeed, the random variables
$S_{k+\ell}+S_{k-\ell}-2S_k$ are symmetric.

The following result is analogous to Theorem~\ref{thm:1a}.

\begin{theorem}\label{thm:4}
Assume that $(X_i)$ is an
i.i.d. sequence\ of
real-valued random variables with distribution\ $F$ {which is regularly
varying} with index
$\alpha>0$, in the sense of (\ref{eq:tb}). Then, with $(a_n)$
defined in (\ref{def:1}), for $f\in\mathcal{F}_\gamma$, {$\gamma
>\max(0,
0.5-\alpha^{-1})$},
\[
\lim_{n\to\infty}P\bigl(a_n^{-1}\widehat
M_n^{(\gamma)}\le x\bigr)= \Phi_\alpha^2(x) ,\qquad x>0 .
\]
\end{theorem}

%s3 ###
\section{Proofs}\label{sec:3}
%s3.1 ###
\subsection[Proof of Lemma 2.4(1)]{Proof of Lemma \protect\ref{lem:2x}(1)}\label{subsecL1}

The following analog of Davis and Resnick \cite{davisresnick1985},
Theorem~2.2, in the case $\mathcal {B}=\mathbb{R}$ is the key to this
result.

\begin{lemma}\label{lem:2} Let $(X_i)$ be an i.i.d. sequence\ of random
elements with values in $\mathcal{B}$. Assume that $X$ is regularly
varying with index
$\alpha>0$ and limit measure\ $\mu$. Then, for any $h\ge1$,
\begin{eqnarray*}
&&\widehat{N}_n=\sum_{t=1}^n \varepsilon_{a_n^{-1}(X_{t},\ldots,X_{t+h-1})}\\
&&\quad\stackrel{d}{\longrightarrow}\quad
\widehat{N}=\sum_{i=1}^\infty\varepsilon_{(J_i,0,\ldots,0)}+ \sum
_{i=1}^\infty
\varepsilon_{(0,J_i,0,\ldots,0)}+\cdots+
\sum_{i=1}^\infty\varepsilon_{(0,\ldots,0,J_i)} ,\qquad n\to\infty,
\end{eqnarray*}
where $\varepsilon_x$ is Dirac measure\ at $x$ and $J_1,J_2,\ldots$ are
the points of a Poisson random measure\ with mean measure\ $\mu$ on
$\mathcal{B}_0$ equipped with the Borel $\sigma$-field. Moreover, on
the right-hand side, in the subscripts of the $\varepsilon$'s, there
are vectors of length $h$. Here, convergence\ in\ distribution\ is in
the space of point measures $M_p$ on $ \mathcal{B}_0^h$, equipped with
the vague topology; see \cite{daleyverejones1988}, Section 9.1.
\end{lemma}

We postpone the proof until the end of this subsection.

\begin{remark}\label{rem:1}
According to Daley and Vere-Jones \cite{daleyverejones1988}, Theorem~9.1.VI, $\widehat
N_n\stackrel{d}{\longrightarrow}
\widehat N$ is
equivalent to the convergence\ of the finite-dimensional distributions
\[
(\widehat N_n(B_1),\ldots,\widehat N_n(B_m)) \stackrel{d}{\longrightarrow
}(\widehat N(B_1),\ldots, \widehat
N(B_m)),\qquad n\to\infty,
\]
for any choice of bounded continuity sets $B_i$ of $\mathcal{B}_0$.
Moreover, according to \cite{daleyverejones1988}, Corollary~9.1.VIII,
it suffices that the sets $B_i$ run through any covering semiring of
bounded continuity {sets} for the limiting measure\ ${ P_{\widehat
N}(\cdot)}=P(\widehat N\in\cdot )$. This means that every open set in
$\mathcal{B}_0$ can be represented as a finite or countable union of
sets from this semiring. Since $\mathcal{B}_0$ is separable, an
important example of such a semiring is obtained by first taking the
open spheres $S(d_k,r_j)$ with centers at the points $d_k$ of a
countable dense set and radii $r_j$ forming a countable dense set in
$(0,1)$, then forming intersections and finally taking proper
differences; see \cite{daleyverejones1988}, page 617. We will make use
of such a semiring in the proof of Lemma \ref{lem:2}.
\end{remark}

A combination of this lemma and the continuous mapping argument
analogous to the one in the proof of Davis and
Resnick \cite{davisresnick1985}, Theorem~2.4, yields
\begin{eqnarray}\label{eq:09}
&&N_n=\sum_{t=1}^n\varepsilon_{a_n^{-1} (X_t,X_t+X_{t+1},\ldots,
X_t+\cdots
+X_{t+h-1})}\nonumber\\[-8pt]\\[-8pt]
&&\quad\stackrel{d}{\longrightarrow}\quad \sum_{i=1}^\infty\varepsilon
_{(J_i,\ldots,J_i)}+
\sum_{i=1}^\infty\varepsilon_{(0,J_i,\ldots,J_i)}+\cdots+
\sum_{i=1}^\infty\varepsilon_{(0,\ldots,0,J_i)}=N ,\qquad n\to\infty.\qquad\nonumber
\end{eqnarray}
Here, the vectors in the subscripts of the $\varepsilon$'s have length $h$.
Write
%
%e3.1 ###
\begin{eqnarray}\label{eq:by}
B(y)=\{(x_1,\ldots,x_h)\in\mathcal {B}^h\dvtx \|x_i\|\le
y ,i=1,\ldots,h\}
\end{eqnarray}
and, for $f\in\mathcal{F}_\gamma$, $\gamma\ge0$,
\[
\widetilde M_{n\ell}^{(\gamma)}=(f(\ell))^{-1}\max_{0\le k \le
n}\|S_{k+\ell}-S_k\| ,\qquad \ell=1,2,\ldots.
\]
Then, for $y>0$, by (\ref{eq:09}), as $n\to\infty$,
%
%e3.2 ###
\begin{eqnarray}\label{oo}
P\bigl(N_n(B(y)^c)=0\bigr)&=& P\bigl(a_n^{-1}\widetilde M_{n1}^{(0)}\le y ,\ldots,
a_n^{-1}\widetilde M_{nh}^{(0)}\le y\bigr)\nonumber\\
&\to& P\bigl(N(B(y)^c)=0\bigr)\\
&=&P \Bigl(\sup_{i\ge1}\|J_i\|\le y ,\sup_{i\ge
1}\|J_i\|\le y ,\ldots,\sup_{i\ge
1}\|J_i\|\le y \Bigr) .\nonumber
\end{eqnarray}
Since $(J_i)$ constitute a Poisson random measure\ on $\mathcal{B}_0$ with
mean measure\ $\mu$, the transformed points $(\|J_i\|)$ constitute a
Poisson random measure\ on $(0,\infty)$ with mean measure\ $\nu$
given by
\[
\nu(y,\infty)=\mu(\{x\in\mathcal{B}_0\dvtx \|x\|>y\})=y^{-\alpha}
\mu(\{x\in\mathcal{B}_0\dvtx \|x\|>1\})=y^{-\alpha} ,\qquad y>0 .
\]
In the last step, we used the definition of $(a_n)$. Moreover, we
assumed that $P(N(\partial B(y)^c)=\break 0)=1$. However,
\begin{eqnarray*}
N(\partial B(y)^c)&=&N(\{x\in\mathcal{B}^h\dvtx \|x_i\|=y ,\|x_j\|\le
y ,j\ne i , \mbox{ for any $i=1,\ldots,h$} \})\\
&\le&\sum_{i=1}^hN(\{x\in\mathcal{B}^h\dvtx \|x_i\|=y\})=0 \qquad\mathrm{a.s.}
\end{eqnarray*}
since the expectation of the right-hand expression is zero.
Hence, $\nu(y,\infty)=y^{-\alpha}$, $y>0$. Writing the points
$\|J_i\|$ in descending order, they have the representation\
\[
\Gamma_1^{-1/\alpha}>\Gamma_2^{-1/\alpha}>\cdots,
\]
where
$(\Gamma_i)$ are the arrivals of a unit-rate homogeneous Poisson
process on
$(0,\infty)$. Therefore, and by (\ref{oo}), we conclude that for
$h\ge
1$, as $n\to\infty$,
\begin{eqnarray*}
&&P\Bigl(a_n^{-1} \max_{\ell\le h} \widetilde M_{n\ell
}^{(\gamma)}\le x\Bigr)\\
&&\quad=
P\Bigl(a_n^{-1} \max_{\ell\le h} (f(\ell))^{-1}\widetilde M_{n\ell
}^{(0)}\le
x\Bigr)\\
&&\quad\to P\Bigl(\sup_{i\ge1} \Gamma_i^{-1/\alpha}\le
x ,(f(2))^{-1}\sup_{i\ge1} \Gamma_i^{-1/\alpha}\le
x ,\ldots,(f(h))^{-1} \sup_{i\ge1} \Gamma_i^{-1/\alpha}\le
x\Bigr)\\
&&\quad=P(\Gamma_1^{-1/\alpha}\le
x)\\
&&\quad=\mathrm{e}^{-x^{-\alpha}}=\Phi_\alpha(x) ,\qquad x>0 .
\end{eqnarray*}
This concludes the proof of Lemma \ref{lem:2x}(1).

\begin{pf*}{Proof of Lemma \ref{lem:2}}
We follow the proofs of Davis and Resnick
\cite{davisresnick1985}, Proposition 2.1 and Theorem~2.2, for the case
$\mathcal{B}=\mathbb{R}$. For $h=1$, the
points of $\widehat N_n$ are independent and\vspace*{-2pt} therefore the convergence
of the
finite-dimensional distributions of $\widehat N_n$ to those of
$\widehat N$ follows from
$\mu_n\stackrel{\widehat w}{\to} \mu$ as $n\to\infty$. Therefore,
we consider
the case $h>1$. We write
\begin{eqnarray*}\widetilde I_n=
\sum_{t=1}^n\varepsilon_{a_n^{-1}(X_t,0,\ldots,0)}+\cdots+
\sum_{t=1}^n\varepsilon_{a_n^{-1}(0,\ldots,0, X_t)} ,
\end{eqnarray*}
where the
points of this process are in $\mathcal{B}_0^h$. Our first aim is to show
that $\widehat N_n(B)-\widetilde I_n(B)\stackrel{P}{\longrightarrow
}0$ as $n\to\infty$ for bounded Borel
sets $ B\subset\mathcal{B}_0^h$ which are bounded away from zero. For\vspace*{-3pt}
this reason (see Remark \ref{rem:1}), it suffices to show that
$\widehat
N_n(B)-\widetilde I_n(B)\stackrel{P}{\longrightarrow}0$ as $n\to
\infty$ for sets $B$ from a covering
semiring of $\mathcal{B}_0^h$. Therefore, it suffices to consider sets
$B=B_1\times\cdots\times B_h$, where each of the sets $B_i\in
\mathcal{B}$ is an element of the semiring generated by the open spheres
$S(d_k,r_j)$, as explained in Remark \ref{rem:1}. We also assume
that $\mu(\partial B_i)=0$, which is possible because $(d_k)$ is
dense in $\mathcal{B}$ and $(r_j)$ in $(0,1)$, and there exist only
{countably} many atoms of $\mu$ because it is finite on bounded
sets.

Since $B$ is bounded away from zero, exactly one of the following
two distinct situations may occur: $(C_1)$: $B$ has no
intersection with the sets $M_i=\{(0,\ldots,0,y,0,\ldots,0)\in
\mathcal{B}^h\dvtx y\in\mathcal{B}\}$, that is, this set consists of the vectors
$(0,\ldots,0,y,0,\ldots,0)$ with $y$ at the $i$th position;
$(C_2)$: $B\cap M_i=B_{i'}$ for $i=i'$ and $B\cap M_i=\varnothing$
for $i\ne i'$. This means that $B$ is either bounded away from the
`axes' $M_i$ or exactly one set $B_i\subset\mathcal{B}$ contains zero.\vspace*{-2pt}
We can now essentially follow the lines of the proof of
Davis and Resnick \cite{davisresnick1985}, Proposition 2.1, in
order to prove $\widehat N_n(B)-\widetilde I_n(B)\stackrel
{P}{\longrightarrow}0$ as $n\to\infty$. (They
prove the result for $\mathcal{B}=\mathbb{R}$ for the more
complicated point processes
involving the points $(t/n,a_n^{-1}(X_t,\ldots,X_{t+h-1})$,
$t=1,2,\ldots.$)

The proof of \cite{davisresnick1985}, Theorem~2.2, can be
adapted by replacing the semiring $S$ in
\cite{davisresnick1985} by the semiring of the sets described
above. Moreover, Davis and Resnick \cite{davisresnick1985}
apply Kallenberg \cite{kallenberg1976}, Theorem~4.2 (which
applies to the convergence\ of point processes on locally compact
spaces). In
our situation, this result can be replaced by the results in
\cite{daleyverejones1988} on the convergence\ of
point processes in a separable complete metric space which were mentioned\vspace*{-2pt}
in Remark \ref{rem:1} above. An adaptation of \cite{davisresnick1985}, Theorem~2.2, yields that $\widehat N_n\stackrel
{d}{\longrightarrow}\widehat N$ as $n\to\infty$.
\end{pf*}

%s3.2 ###
\subsection[Proof of Lemma 2.4(2)]{Proof of Lemma \protect\ref{lem:2x}(2)}\label{subsecL2}
By the definition of the class $\mathcal{F}_\gamma$, it suffices to
prove the result for the functions $f(\ell)=\ell^\gamma$. We will show
that {the quantities $a_n^{-1}\widetilde M_{n\ell}^{(\gamma)}$ for
large values of $\ell$ do not contribute} to the limit distribution\ of
$a_n^{-1} \widetilde M_n^{(\gamma)}$. The key to the proof is the
following inequality.
%{\bf Replace $H<n/2$ by $n$. Possible?}

\begin{lemma}\label{lem:ineq} Let $(X_i)$ be an i.i.d. sequence\ with
values in
$\mathcal{B}$. Then, for any
$\delta,\gamma>0$, $h\ge1$ and $H\le n$,
\[
P \Bigl(\max_{h\le\ell\le H} \ell^{-\gamma}\sup_{k\le
n}\|S_{k+\ell}-S_k\|>\delta a_n \Bigr) \le2 \sum_{j=J_1}^{J_0} 2^j
P \Bigl( \max_{1\le k\le2n2^{-j}}\|S_{k}\|> \delta(n
2^{-j})^\gamma a_n \Bigr) ,
\]
where $J_0=\log_2(n/h)$, $J_1=\log_2(n/H)+1$ and $\log_2 x$
denotes the dyadic logarithm.
\end{lemma}

Here, and in what follows, we abuse notation when we write
$\sum_{j=a}^b x_j$ instead of $\sum_{j:a\le j\le b}x_j$ for real
values $a<b$.

\begin{pf*}{Proof of Lemma \ref{lem:ineq}}
We use a dyadic splitting of the $\ell$- and $k$-index ranges.
Recall the definitions of $J_0$, $J_1$, { where we assume, for
simplicity, that these numbers are integers.} Setting
\[
I_j=(n2^{-j}, n2^{-j+1}] ,\qquad j=J_1,\ldots,J_0 ,
\]
we obtain
\[
\bigcup_{j=J_1}^{J_0}I_j=\{h, h+1, \dots, H\}
\]
and, therefore,
%
%e3.3 ###
\begin{eqnarray}\label{eq:03}
\max_{h\le\ell\le H}
\ell^{-\gamma}\max_{k\le n} \|S_{k+\ell}-S_k\|&=& \max_{J_1\le
j\le J_0} \max_{\ell\in I_j} \ell^{-\gamma}\max_{1\le k\le n}
\|S_{k+\ell}-S_k\|\nonumber\\
&\le& \max_{J_1\le j\le J_0} (n^{-1}2^{j})^\gamma\max_{\ell\in I_j}
\max_{1\le k\le n}\|S_{k+\ell}-S_k\|\\
&\le&\max_{J_1\le j\le J_0}(n^{-1}2^{j})^\gamma\max_{\ell\in
I_j}\max_{1\le i< 2^j}\max_{(i-1)n2^{-j}\le k<
in2^{-j}}\|S_{k+\ell}-S_k\| .\nonumber
\end{eqnarray}
We observe that for $n2^{-j}< \ell\le n2^{-j+1}$ and
$(i-1)n2^{-j}\le k< in2^{-j}$,
\begin{eqnarray*}\|S_{k+\ell}-S_k\| & \le&
\bigl\|S_{k+\ell}-S_{[in2^{-j}]}\bigr\|+
\bigl\|S_{[in2^{-j}]}-S_k\bigr\|\\
& \le&
\max_{in2^{-j}<u< (i+2)n2^{-j}}\bigl\|S_{u}-S_{[in2^{-j}]}\bigr\|+
\max_{(i-1)n2^{-j}\le k< in2^{-j}}\bigl\|S_{[in2^{-j}]}-S_k\bigr\| .
\end{eqnarray*}
Hence, by virtue of (\ref{eq:03}), we obtain the bound
\[
P \Bigl(\max_{h\le\ell\le H}\ell^{-\gamma}\max_{k\le n}
\|S_{k+\ell}-S_k\|>\delta
a_n \Bigr)\le P_1+P_2 ,
\]
where
\begin{eqnarray*}
P_1&=& P \Bigl(\max_{J_1\le j\le J_0} (n^{-1}2^{j})^\gamma\max_{1\le i<
2^j} \max_{in2^{-j}<u<(i+2)n2^{-j}}\bigl\|S_{u}-S_{[in2^{-j}]}\bigr\|>\delta
a_n\Bigr) ,\\
P_2&=&P \Bigl(\max_{J_1\le j\le J_0} (n^{-1}2^{j})^\gamma
\max_{1\le i< 2^j} \max_{(i-1)n2^{-j}\le k<
in2^{-j}}\bigl\|S_{[in2^{-j}]}-S_k\bigr\|>\delta a_n \Bigr) .
\end{eqnarray*}
Finally, we obtain
\begin{eqnarray*}
P_1&\le& \sum_{j=J_1}^{J_0} \sum_{1\le i<
2^j} P \Bigl(\max_{in2^{-j}<u<(i+2)n2^{-j}}\bigl\|S_{u}-S_{[in2^{-j}]}\bigr\|> \delta
(n2^{-j+2})^\gamma a _n \Bigr)\\
&=&\sum_{j=J_1}^{J_0} 2^j P \Bigl( \max_{1\le k\le
2n2^{-j}}\|S_{k}\|> \delta(n2^{-j})^\gamma a_n \Bigr) .
\end{eqnarray*}
In the last step, we used the i.i.d. property of $(X_i)$. The
corresponding bound for $P_2$ is similar.\quad\mbox{}
\end{pf*}

We are now ready for the second part of Lemma \ref{lem:2x}. We
consider the truncated elements
\[
X_i'=X_iI_{\{\|X_i\|\le h^\gamma a_n\}} ,\qquad
\widetilde{X}_i=X_i'-EX_i' ,\qquad i=1,\ldots,n ,
\]
and the corresponding partial sums $S_k'=\sum_{i=1}^k X_i'$ and
$\widetilde S_k=\sum_{i=1}^k\widetilde X_i,$ $k=1, \dots, n$, with
$S_0'=\widetilde
S_0=0.$ By virtue of Lemma \ref{lem:ineq}, we conclude that for
any $\delta>0$,
%
%e3.4 ###
\begin{eqnarray}\label{prop1:1}
&&P \Bigl(\max_{h\le\ell\le n}\ell^{-\gamma}\max_{0\le
k\le
n-\ell}\|S_{k+\ell}-S_k\|> \delta a_n \Bigr)\nonumber\\
&&\quad\le P\Bigl(\max_{1\le k\le
n}\|X_k\|\ge h^\gamma a_n\Bigr)
+P \Bigl(\max_{h\le\ell\le n}\ell^{-\gamma}\max_{0\le k\le
n-\ell}\|S'_{k+\ell}-S'_k\|> \delta
a_n \Bigr)\\
&&\quad\le P\Bigl(\max_{1\le k\le n}\|X_k\|\ge h^\gamma
a_n\Bigr)+2\sum_{j=1}^{\log_2(n/h)} 2^j { Q_j ,}\nonumber
\end{eqnarray}
where
\[
 Q_j =P \Bigl( \max_{1\le k\le2n2^{-j}}\|S'_{k}\|> \delta
(n2^{-j})^\gamma a_n \Bigr) .
\]
Since (\cite{gnedenko1943}; see, for example,
\cite{embrechtskluppelbergmikosch1997}, Theorem~3.3.7,
for a more recent reference)
\[
\lim_{h\to\infty}\lim_{n\to\infty}P\Bigl(\max_{1\le k\le n}\|X_k\|\ge
h^\gamma a_n\Bigr)=1-\lim_{h\to\infty}\mathrm{e}^{-h^{-\gamma\alpha
}}=0 ,
\]
it suffices to show that
%
%e3.5 ###
\begin{equation}\label{prop1:3}
\lim_{h\to\infty}\limsup_{n\to\infty}\sum_{j=1}^{\log_2(n/h)} 2^j
Q_j=0 .
\end{equation}
Write
\[
\Delta_{nj}=a_n(n2^{-j})^\gamma\quad\mbox{and}\quad
N=[2n2^{-j}] .
\]

First, we consider the case $\alpha>1.$ By assumption, $EX=0$ and,
therefore, we have
\[
\max_{1\le k\le N}\|ES_k'\| =N\bigl\|EX I_{\{\|X\|\ge
h^\gamma a_n\}}\bigr\| \le N E\bigl(\|X\|I_{\{\|X\|\ge h^\gamma
a_n\}}\bigr) .
\]
Since $\|X\|$ is regularly varying\ with index $\alpha$,
an application of Karamata's theorem yields that as $n\to\infty$,
%
%e3.6 ###
\begin{equation}\label{prop1:4}
E\bigl(\|X\|I_{\{\|X\|\ge h^\gamma a_n\}}\bigr)\sim{ c_\alpha} n^{-1}a_n
h^{\gamma(1-\alpha)} .
\end{equation}
Hence, since $N\le n$ and, therefore,
$cn^{-1}N^{1-\gamma}h^{\gamma(1-\alpha)}\le\delta/2$ for large
$n$ {and some constant $c>0$}, we have
\[
Q_j\le P \Bigl( \max_{1\le k\le N}\|S'_{k}-ES_k'\|> \delta N^\gamma
a_n-cNh^{\gamma(1-\alpha)} a_n n^{-1} \Bigr)\le\widetilde Q_j ,
\]
where
\[
\widetilde Q_j=P \Bigl( \max_{1\le k\le N}\|\widetilde S_{k}\|> (\delta/2)
N^\gamma a_n \Bigr) .
\]
Since the sequence\ $(\|\widetilde S_k\|)_{k=0, 1, \dots}$ constitutes a
submartingale, an application of the Chebyshev and Doob
inequalities for $p>1$ yields, {for some
constant $c>0$}, that
%
%e3.7 ###
\begin{equation}\label{prop1:5}
\widetilde Q_j\le(\delta/2)^{-p}\Delta_{nj}^{-p}E \Bigl(\max_{1\le k\le
N}\|
\widetilde
S_k\|^p \Bigr)\le c\Delta_{nj}^{-p}E\|\widetilde S_N\|^p .
\end{equation}
We proceed by applying an $L_p$-inequality for sums of independent
mean zero random elements (see \cite{ledouxtalagrand1993}, Theorem~6.20). We obtain, for $p>2$,
\begin{equation}\label{eq:nk}
E\|\widetilde{S}_N\|^p\le c [ (E\|\widetilde S_N\|)^p+NE\|\widetilde
X_1\|^p]
\end{equation}
with a constant $c$ depending on $p$ only. For fixed $\gamma>0$
and $\alpha>1$, let us choose $\beta>0$ such that $\beta<\alpha$
and $\gamma>\beta^{-1}-\alpha^{-1}.$ We then have
\begin{eqnarray*}
E\|\widetilde S_N\|&=&E \Biggl\|\sum_{i=1}^NX_i'-N EX_1' \Biggr\|=
E \Biggl\|S_N-\sum_{i=1}^NX_iI_{\{\|X_i\|>h^\gamma a_n\}}-
NE\bigl(XI_{\{\|X\|>h^\gamma a_n\}}\bigr) \Biggr\|\nonumber\\
&\le& E\|S_N\|+2N E\bigl(\|X\|I_{\{\|X\|>h^\gamma a_n\}}\bigr) .
\end{eqnarray*}
By (\ref{prop1:4}) and assumption (\ref{restr:1}), we conclude that
%
%e3.8 ###
\begin{equation}\label{eq:98}
E\|\widetilde S_N\|\le c\bigl[N^{1/\beta}+Nn^{-1}a_n h^{\gamma(1-\alpha)}\bigr] .
\end{equation}
Again, by regular variation\ of $\|X\|$ and Karamata's theorem, for
$p>\max(2,\alpha)$,
as $n\to\infty$,
%
%e3.9 ###
\begin{equation}\label{prop1:8}
E\|\widetilde X_1\|^p
%=\int_0^{h^\gamma a_n}t^{p-1}P(\|X_1\|>t)dt=
%a_n^p\int_0^{h^\gamma}t^{p-1}P(\|X_1\|>ta_n)dt
\sim{ c_\alpha} a_n^{p}n^{-1}h^{\gamma(p-\alpha)} .
\end{equation}
Combining (\ref{prop1:5})--(\ref{prop1:8}), we obtain
\begin{eqnarray*}
\sum_{j=1}^{\log_2(n/h)}2^j\widetilde Q_j&\le& c\sum_{j=1}^{\log
_2(n/h)} 2^j
a_n^{-p}N^{-p\gamma} \bigl[N^{p/\beta}+N^pn^{-p}a^p_n
h^{p\gamma(1-\alpha)}+
Na_n^{p}n^{-1}h^{\gamma(p-\alpha)} \bigr]\\
&\le& c [I_1+I_2+I_3] ,
\end{eqnarray*}
where
\begin{eqnarray*}
I_1&=&
a_n^{-p}n^{-p\gamma+p/\beta}\sum_{j=1}^{\log
_2(n/h)}2^{j(1{+}p\gamma{ -}p/\beta)} ,\\
I_2&=&n^{-p\gamma}h^{p(1-\alpha)}\sum_{j=1}^{\log
_2(n/h)}2^{j-pj+p\gamma j} ,\\
I_3&=&n^{-p\gamma} h^{\gamma(p-\alpha)}\sum_{j=1}^{\log
_2(n/h)}2^{p\gamma
j} .
\end{eqnarray*}
If $\gamma\le1/\beta$, then {using the fact that} $p>\max(2,\alpha
)$, for
some constants $c>0$ and a slowly varying\ function\ $\ell$,
\[
I_1\sim{c} a_n^{-p}n^{-p\gamma+p/\beta}(n/h)^{1{ +}p\gamma{ -}p/\beta}=
c h^{-1-p\gamma+p/\beta}(\ell(n))^{-p}n^{-p/\alpha+1}=\mathrm{o}(1),\qquad n\to\infty.
\]
If $\gamma> 1/\beta-1/p$, then $I_1=\mathrm{o}(1)$ as $n\to\infty$ by
choosing $p>1/(\gamma-1/\beta)$. Next, we see that $I_2=\mathrm{O}(n^{-p+1})$
as $n\to\infty$ if $\gamma\ge1$ and $I_2= \mathrm{O}(n^{-p\gamma})$ as
$n\to\infty$ if
$\gamma<1$ and $p>1/(1-\gamma).$ Finally, $I_3\le c
h^{-\gamma\alpha}$ for some $c>0$ and the right-hand side\ converges
to zero
as $h\to\infty.$ This proves (\ref{prop1:3}) for $\alpha>1.$

The case $\alpha=1$,
$EX=0$, can be handled following the lines of the proof above. Then,
(\ref{prop1:4}) does not remain valid. However, Karamata's theorem
yields that $f_1(x)=E\|X\|_{\{\|X\|>x\}}$ is a slowly varying\
function. This
fact suffices to derive the corresponding relations after
(\ref{prop1:4}).

We now consider the cases $0<\alpha<1$ and $\alpha=1$,
$E\|X\|=\infty$. %We start with the case $\alpha<1$ and comment on
%the case $\alpha=1$ later.}
We have
\[
\sum_{j=1}^{\log_2(n/h)} 2^j Q_j\le\sum
_{j=1}^{\log_2(n/h)}
2^jP (T_{N}-N E\|X'\|>\delta N^\gamma a_n- N E\|X'\| ) ,
\]
where $T_k=\sum_{i=1}^k \|X_i'\|$. Another application of
Karamata's theorem yields, as $n\to\infty$,
\[
E\bigl(\|X\|I_{\{\|X\|\le h^\gamma a_n\}}\bigr)
%&=\int_0^{h^\gamma
%&a_n}P(\|X_1\|>t)dt= a_n\int_0^{h}P(\|X_1\|>ta_n)dt\\
\sim\
\cases{  c_\alpha a_n n^{-1}h^{\gamma(1-\alpha)},
&\quad when $\alpha\not=1$,\cr
\mbox{slowly varying,}&\quad when
$\alpha=1$.
}
\]
We now easily deduce that $NE\|X'\|=\mathrm{o}(N^\gamma a_n)$ as $n\to\infty$ for
$\gamma>0$. Thus, we have, for large~$n$, by Kolmogorov's inequality,
\begin{eqnarray*}
P (T_{N}-N E\|X'\|>\delta N^\gamma a_n- N E\|X'\| )&\le&
P \bigl(T_{N}-N E\|X'\|>(\delta/2)
N^\gamma a_n \bigr)\\
&\le &c_\delta N^{-2\gamma}a_n^{-2}N\operatorname{var}(\|X'\|) .
\end{eqnarray*}
%
%Then we have
%$$
%2^j N E(\|X\|I_{\{\|X\|\le h^\gamma a_n\}}).
%$$
%Another application of Karamata's theorem yields as $\nto$
%E(\|X\|I_{\{\|X\|\le h^\gamma a_n\}})
%%&=\int_0^{h^\gamma
%%&a_n}P(\|X_1\|>t)dt= a_n\int_0^{h}P(\|X_1\|>ta_n)dt\\
By Karamata's theorem, $\operatorname{var}(\|X'\|)\sim c_\alpha a_n^2 n^{-1}
h^\gamma$ as $n\to\infty$. So, for large $n$,
\[
\sum_{j=1}^{\log_2(n/h)} 2^j Q_j\sim c_{\alpha, \delta}
\sum_{j=1}^{\log_2(n/h)} 2^jN^{1-2\gamma}n^{-1}h^\gamma\le
c_{\alpha, \delta} \sum_{j=1}^{\log_2(n/h)} 2^{2\gamma j}
n^{-2\gamma}h^\gamma\le c_{\alpha, \delta} h^{-\gamma}
\]
and we conclude that (\ref{prop1:3}) indeed holds.
%and some constant $c>0$,
%cn^{-1}h^{1-\alpha}\sum_{j=1}^{\log_2(n/h)}2^jN^{1-\gamma}
%j}\le ch^{-\gamma\alpha} .
%Since the \rhs\ converges to zero as $h\to\infty$ the proof of
%(\ref{prop1:3}) is complete {\green for $\alpha<1$.
%For $\alpha=1$, $E\|X\|=\infty$, Karamata's theorem yields that
%$f_2(x)=E(\|X\|I_{\{\|X\|\le x\}})$ is a \slvary\ \fct\ and therefore
%the argument above does not apply. However, applying the triangle
%inquality, we obtain
%N^\gamma a_n- N E\|X'\| ) ,
%where $T_k=\sum_{i=1}^k \|X_i'\|$. Since $E\|X'\|=f_2(h^\gamma a_n)$ is
%a \slvary\ \fct\ of $n$ and $N E\|X'\|\le n E\|X'\|=o(N^\gamma
%a_n)$ as $\nto$ for $\gamma>0$, we have $\delta N^\gamma a_n-
%N E\|X'\|\ge(\delta/2) N^\gamma a_n$ for large $n$.
%An application of Prokohorov's inequality (e.g. Ledoux and Talagrand
%&\le& \exp\{- (\delta/2) (N/h)^\gamma(\log(1+ (\delta/2)(hN)^
%By Karamata's theorem, $\var(\|X'\|)\sim c_\alpha a_n^2 n^{-1} h^
%$\nto$.
%Therefore and since
%the \rhs\ is bounded by $(n N^{\gamma-1})^{-c(\delta)(N/h)^\gamma}$
%for some constant $c(\delta)>0$. We conclude that
%$\sum_{j=1}^{\log_2(n/h)} 2^j Q_j \to0$ as $\nto$.
%
%}
This completes the proof of the lemma.
% \hfill\halmos

%s3.3 ###
\subsection[Proof of Theorem 2.2]{Proof of Theorem~\protect\ref{thm:1}}
The proof of (\ref{main:1}) is immediate from Lemma
\ref{lem:2x}. It thus {suffices to prove} (\ref{main:2}). We
achieve this by showing that the sequences $(a_n^{-1}\widetilde
M_n^{(\gamma)})$ and $(a_n^{-1}\widetilde T_n^{(\gamma)})$ have the same
asymptotic\ behavior.

Throughout the proof, we set
\[
V_\ell(i, j)=\max_{i\le k\le
j}\|S_{k+\ell}-S_k-\ell\overline{X}_n\| ,\qquad 0\le i<j\le n .
\]
The argument of Remark \ref{rem:xu} allows us to assume that $\gamma
\in(0, 1)$ and we start
by observing that, in view of Remark \ref{rem:97}, the sequences
$(a_n^{-1}\widetilde M_n^{(\gamma)})$ and $(a_n^{-1}\zeta_n^{(\gamma)})$
have the same limit distribution. Next, we observe that
$(a_n^{-1}\zeta_n^{(\gamma)})$ has the same limit distribution as
\[
a_n^{-1}\zeta_n'^{(\gamma)}= a_n^{-1}\max_{1\le\ell\le d_n}
\bigl(f\bigl(\ell(1-\ell/n)\bigr)\bigr)^{-1} V_\ell(0, n-\ell) ,\qquad % \max_{0\le k\le
%n-\ell}\|S_{k+\ell}-S_k-\ell\ol{X}_n\| ,
n\ge1 ,
\]
for any sequence $d^2_n\to\infty$ such that $d_n/n\to0$ as
$n\to\infty$. Indeed, we have
\[
\inf_{1\le\ell\le d_n}\bigl(f\bigl(\ell(1-d_n/n)\bigr)/f(\ell)\bigr) \zeta_n'^{(\gamma)}\le
\zeta_n^{(\gamma)}\le\max\bigl(\zeta_n'^{(\gamma)}, \Delta_n\bigr) ,
\]
where
$\Delta_n=\max_{\ell\ge d_n}(f(\ell))^{-1}V_\ell(0,
n-\ell)=\mathrm{o}_P(a_n)$ due to Lemma \ref{lem:2x} and Remark
\ref{rem:97}. {By the definition of the class $\mathcal{F}_\gamma$, we
have that $\inf_{1\le\ell\le d_n}(f(\ell(1-d_n/n))/f(\ell))\to1$ as
$n\to\infty.$} Again by Lemma \ref{lem:2x} and Remark
\ref{rem:97}, we conclude that the sequence
$(a_n^{-1}\zeta_n'^{(\gamma)})$ has the same asymptotic\ distribution as
$(a_n^{-1}\zeta_n''^{(\gamma)})$, where
\[
\zeta_n''^{(\gamma)}= \max_{1\le\ell\le n/2}
\bigl(f\bigl(\ell(1-\ell/n)\bigr)\bigr)^{-1} V_\ell(0, n-\ell) ,\qquad %\max_{0\le k\le
%n-\ell}\|S_{k+\ell}-S_k-\ell\ol{X}_n\| .
n\ge1 .
\]
Finally, we show that
%
%e3.10 ###
\begin{equation}\label{last:1}
\Delta_n=a_n^{-1}\max_{n/2<\ell< n} \bigl(\ell(1-\ell/n)\bigr)^{-\gamma}
V_\ell(0, n-\ell)\stackrel{P}{\longrightarrow}
0 \qquad\mbox{as } n\to\infty.
\end{equation}
We use the following identity:
\begin{eqnarray*}
S_{k+\ell}-S_k-\ell\overline{X}_n&=&\sum_{i=k+1}^{k+\ell
}(X_i-\overline{X}_n)\\
&=&
- \Biggl[\sum_{i=k+\ell+1}^n(X_i-\overline{X}_n)+\sum_{i=1}^k
(X_i-\overline{X}_n) \Biggr].
\end{eqnarray*}
In view of the identical distributions of the $X_i$'s, the proof of
(\ref{last:1}) reduces to showing that, as $n\to\infty$,
%
%e3.11 ###
\begin{equation}\label{last:2}
\Delta_n'=a_n^{-1}\max_{n/2<\ell< n} \bigl(\ell(1-\ell/n)\bigr)^{-\gamma}
\max_{0\le k\le
n-\ell} \Biggl\|\sum_{i=k+1}^{k+n-\ell}(X_i-\overline{X}_n) \Biggr\| \stackrel
{P}{\longrightarrow}
0 .
\end{equation}
By virtue of (\ref{eq:67}), we have
\begin{eqnarray*}
\Delta_n'&=&a_n^{-1}\max_{1\le\ell<n/2}\bigl(\ell(1-\ell/n)\bigr)^{-\gamma}
\max_{0\le k\le
\ell} \Biggl\|\sum_{i=k+1}^{k+\ell}(X_i-\overline{X}_n) \Biggr\|\\
&\le&2a_n^{-1}\max_{1\le\ell<n/2}\bigl(\ell(1-\ell/n)\bigr)^{-\gamma}
\max_{0\le k\le2\ell} [\|S_k\|+\ell\|\overline{X}_n\| ]\\
&\le&2^{\gamma+1}a_n^{-1}\max_{1\le
\ell<n/2}\ell^{-\gamma}\max_{0\le k\le2\ell}\|S_k\|+\mathrm{o}_P(1)\\
&\le&2^{2\gamma+1}a_n^{-1}\max_{1\le k\le
n}\|k^{-\gamma}S_k\|+\mathrm{o}_P(1) ,\qquad n\to\infty.
%2^{2\gamma+2}a_n^{-1}\max_{1\le k\le n}\|k^{-\gamma}S_k\|.
\end{eqnarray*}
By assumption (\ref{restr:1}), choosing $\beta=2$ if
$\gamma>0.5-1/\alpha$ or $\beta<\alpha$ such that
$\gamma>1/\beta-1/\alpha$, we have
\begin{eqnarray*}
a_n^{-1} \max_{1\le k\le n}k^{-\gamma}\|S_k\|&=&
a_n^{-1}\max_{1\le k\le
n}k^{-\gamma+1/\beta}\|k^{-1/\beta}S_k\|\\& \le&
a_n^{-1} \max(1,n^{-\gamma+1/\beta})\max_{1\le k\le
n}\|k^{-1/\beta}S_k\|\stackrel{P}{\longrightarrow}0 ,\qquad n\to\infty.
\end{eqnarray*}
This concludes the proof of the theorem.
\subsection[Proof of Theorem 2.10]{Proof of Theorem~\protect\ref{thm:1a}}
The proof is similar to that of Theorem~\ref{thm:1}.
Lemma \ref{lem:2} remains valid with $a_n$ replaced by
$b_n=p^{1/\alpha} a_n$, but the limiting Poisson random measure\ with
state space $\mathbb{R}\backslash\{0\}$ has mean measure~$\mu$
given by
$\mu(x,\infty)=x^{-\alpha}$ and $\mu(-\infty,-x)=(q/p)x^{-\alpha}$
for $x>0$.  Consider the set
\[
B(x,y)^c= (-\infty,-x)\cup
(y,\infty) ,\qquad x,y>0 .
\]
Recall the definition of $N_n$
(with $a_n$ replaced by $b_n$) from (\ref{eq:09}). Then, using
Lemma \ref{lem:2} and the same ideas as in the proof of
Lemma \ref{lem:2x}(1), for $h\ge1$ and $x,y>0$,
\begin{eqnarray*}
&&P\bigl(N_n(B(x,y)^c)=0\bigr)\\
&&\quad= P \Bigl(b_n^{-1}\max_{\ell=1,\ldots,h}(f(\ell))^{-1} \max_{k=0,\ldots,n-\ell}
(S_{k+\ell}-S_k)\le y ,\\
&&\qquad\hphantom{P \Bigl(} b_n^{-1}\min_{k=1,\ldots,h}(f(\ell))^{-1}
\min_{k=0,\ldots,n-\ell} (S_{k+\ell}-S_k)\ge-x \Bigr)\\
&&\quad\to \exp\bigl\{-\mu\bigl((-\infty,-x)\cup
(y,\infty) \bigr) \bigr\}\\
&&\quad=\exp\{-(q/p) x^{-\alpha}- y^{-\alpha}\}\\
&&\quad=\Phi_\alpha^{q/p}(x)\Phi_\alpha(y) .
\end{eqnarray*}
Furthermore, for fixed $h\ge1$,
\begin{eqnarray*}
&&P \Bigl(b_n^{-1}\max_{\ell=1,\ldots,h}(f(\ell))^{-1}
\max_{k=0,\ldots,n-\ell} (S_{k+\ell}-S_k)\le y ,\\
&&\hphantom{P \Bigl(}b_n^{-1}
\min_{k=1,\ldots,h}(f(\ell))^{-1}
\min_{k=0,\ldots,n-\ell} (S_{k+\ell}-S_k)\le-x \Bigr)\\
&&\quad=P\Bigl(b_n^{-1}\max_{\ell=1,\ldots,h}(f(\ell))^{-1}
\max_{k=0,\ldots,n-\ell} (S_{k+\ell}-S_k)\le y\Bigr)\\
&&\qquad{}-P\Bigl(b_n^{-1}\max_{\ell=1,\ldots,h}(f(\ell))^{-1}
\max_{k=0,\ldots,n-\ell} (S_{k+\ell}-S_k)\le y ,\\
&&\qquad\hphantom{{}-P\Bigl(} b_n^{-1}
\min_{k=1,\ldots,h}(f(\ell))^{-1}
\min_{k=0,\ldots,n-\ell} (S_{k+\ell}-S_k)> -x\Bigr)\\
&&\quad\to \Phi_\alpha(y)\bigl(1- \Phi_\alpha^{q/p}(x)\bigr) ,\qquad x,y>0 .
\end{eqnarray*}
The right-hand side\ can be extended to a {bivariate} distribution\ in {a}
natural way. An application of Lemma \ref{lem:2x}(2) shows that
this distribution\ is the joint limit distribution\ of
$b_n^{-1}(m_n^{(\gamma)},M_n^{(\gamma)})$.
% \hfill\halmos

%s3.5 ###
\subsection[Proof of Theorem 2.12]{Proof of Theorem~\protect\ref{thm:1b}} One can follow
the lines of the proof of (\ref{main:2}) {to show} that the
sequences $(b_n^{-1}M_n^{(\gamma)})$ and
$(b_n^{-1}T_n^{(\gamma)})$ have the same limiting
distribution.
% \hfill\halmos

%s3.6 ###
\subsection[Proof of Theorem 2.13]{Proof of Theorem~\protect\ref{thm:4}}
The proof is similar to that of Theorem~\ref{thm:1a}. We sketch the
main ideas.
We first observe that the symmetric random variable\ $\widehat
X=X_1-X_2$ is regularly varying:
\[
P(X_1-X_2>x)\sim P(X>x)+P(X< -x)= P(|X|>x) ,\qquad x\to\infty;
\]
see \cite{embrechtskluppelbergmikosch1997},
Lemma A.3.26. Hence, $n P(X_1-X_2>a_n)\to1$ as $n\to\infty$.

Following the lines of the proof of Theorem~\ref{thm:1},
one can show that
\begin{eqnarray*}
\lim_{h\to\infty}\limsup_{n\to\infty} P \Bigl(a_n^{-1}\max_{\ell
=h,\ldots,n}(f(\ell))^{-1}
\max_{k=\ell+1,\ldots,n-\ell} |S_{k+\ell}+S_{k-\ell}-2S_k|>\delta
\Bigr)=0 ,\qquad \delta>0 .
\end{eqnarray*}
Hence, it suffices to show that for any fixed $h\ge1$,
\begin{eqnarray*}
\lim_{n\to\infty}P \Bigl(a_n^{-1}\max_{\ell=1,\ldots,h}
(f(\ell))^{-1}\max_{k=\ell+1,\ldots,n-\ell}
(S_{k+\ell}+S_{k-\ell}-2 S_k)\le x \Bigr)=\Phi_\alpha^2(x) ,\qquad
x>0 .
\end{eqnarray*}
This is again achieved by a point  process argument in the
spirit of
Davis and Resnick~\cite{davisresnick1985}. The same argument as in
Section \ref{subsecL1} yields, for any fixed $\ell\ge1$,
\begin{eqnarray*}
&&\sum_{k=\ell+1}^n\varepsilon_{a_n^{-1}(X_{k+1}-X_{k-1},(X_{k+1}-X_{k-1})+
(X_{k+2}-X_{k-2}),\ldots,(X_{k+1}-X_{k-1})+\cdots+(X_{k+\ell
}-X_{k-\ell}))}\\
&&\quad\stackrel{d}{\longrightarrow} \sum_{i=1}^\infty
\bigl[\varepsilon_{(J_i,,\ldots,J_i)}+\varepsilon_{(-J_i,\ldots,-J_i)}
\bigr]+\sum_{i=1}^\infty
\bigl[\varepsilon_{(0,J_i,\ldots,J_i)}+\varepsilon_{(0,-J_i,\ldots,-J_i)}
\bigr]+\cdots
\\
&&\quad\hphantom{\stackrel{d}{\longrightarrow}}{}+\sum_{i=1}^\infty
\bigl[\varepsilon_{(0,\ldots,0,J_i)}+\varepsilon_{(0,\ldots,0,-J_i)} \bigr]
,\qquad n\to\infty,
\end{eqnarray*}
where $(J_i)$ are the points of a Poisson random measure\ on $\mathcal{B}_0$
with mean measure\ $\mu$ satisfying $\mu(x,\infty)=\mu(-\infty
,-x]=x^{-\alpha}$, $x>0$.
This limit result implies that for $h\ge1$ and $x>0$,
\begin{eqnarray*}
&&P\Bigl(a_n^{-1}\max_{l=1,\ldots,h} (f(\ell))^{-1}\max
_{k=l+1,\ldots,n-\ell}
(S_{k+l}+S_{k-l}-2S_k)
\le x\Bigr)\\
&&\quad\to P\Bigl(\sup_{i\ge1 }|J_i|\le x,(f(2))^{-1}\sup_{i\ge1}|J_i|\le
x,\ldots,(f(h))^{-1}\sup_{i\ge1}|J_i|\le x\Bigr)\\
&&\quad=\Phi_\alpha^2(x) ,\qquad n\to\infty.
\end{eqnarray*}
This concludes the proof.
% \hfill\halmos

\section*{Acknowledgements}
We would like to thank Herold
Dehling for numerous discussions on the proofs of the results and
their presentation. We would like to thank one of the referees
for a very careful report which led to a substantial
improvement of the presentation of this paper.
The research of Thomas Mikosch was supported in part by the Danish
Research Council (FNU) Grant 272-06-0442.

\printhistory

\end{document}